\documentclass[review]{elsarticle}
\usepackage{amsmath,amsfonts,amssymb,amsthm,bm,accents,statex}
\usepackage{lineno}
\usepackage[margin=1in]{geometry}
\usepackage{hyperref}
\hypersetup{
	colorlinks=false,
	pdfborder={0 0 0},
}
\modulolinenumbers[5]

\journal{Statistics and Probability Letters}

\newtheorem{theorem}{Theorem}
\newtheorem{assumption}{Assumption}

\newtheorem{lemma}{Lemma}
\newtheorem{proposition}{Proposition}
\newtheorem{remark}{Remark}

\biboptions{authoryear}

\begin{document}

\begin{frontmatter}

\title{Limit Theory for Moderate Deviation from Integrated GARCH Processes}

\author{Yubo Tao\fnref{myft1}}
\address{90 Stamford Rd, Singapore Management University}
\fntext[myft1]{I would like to thank the co-editor, an associate editor and the referees for helping improve the paper. All possible errors are mine. Yubo Tao, School of Economics, Singapore Management University, 90 Stamford Road, Singapore 178903. Email: yubo.tao.2014@phdecons.smu.edu.sg.}


\begin{abstract}
	This paper develops the limit theory of the GARCH(1,1) process that moderately deviates from IGARCH process towards both stationary and explosive regimes. The GARCH(1,1) process is defined by equations $u_t = \sigma_t \varepsilon_t$, $\sigma_t^2 = \omega + \alpha_n u_{t-1}^2 + \beta_n\sigma_{t-1}^2$ and $\alpha_n + \beta_n$ approaches to unity as sample size goes to infinity. The asymptotic theory extends \cite{BerkesHorvathKokoszka2005} by allowing the parameters to have a slower rate of convergence. The results can be applied to unit root test for processes with mildly-integrated GARCH innovations (e.g. \cite{Boswijk2001}, \cite{Cavaliere2007, Cavaliere2009}) and deriving limit theory of estimators for models involving mildly-integrated GARCH processes (e.g. \cite{JensenRahbek2004}, \cite{FrancqZakoian2012,FrancqZakoian2013}).
\end{abstract}


\end{frontmatter}


\section{Introduction} \label{secIntro}

The model considered in this paper is a GARCH(1,1) process:
\begin{align*}
&\text{(Return Process)} \quad && u_t = \sigma_t \varepsilon_t, \\
&\text{(Volatility Process)} \quad && \sigma_t^2 = \omega + \alpha_n u_{t-1}^2 + \beta_n\sigma_{t-1}^2,\quad
\omega > 0\text{, }\alpha_n \geq 0\text{, and }\beta_n \geq 0,
\end{align*}
where $\lbrace\varepsilon_t\rbrace_{t=0}^{n}$ is a sequence of independent identically distributed (i.i.d) variables such that $E\varepsilon_0 = 0$ and $E\varepsilon_0^2 = 1$. 

Unlike conventional GARCH(1,1) process, the innovation process considered in this paper is a mildly-integrated GARCH process whose key parameters, $\alpha_n$ and $\beta_n$, are changing with the sample size, viz.
\begin{equation*}
	\alpha_n = O(n^{-p}), \quad \beta_n = 1+ O(n^{-q}), \text{ where } p, q \in (0, 1),
\end{equation*}
and 
\begin{equation*}
	\gamma_n = \alpha_n + \beta_n - 1 = O(n^{-\kappa}), \quad \kappa = \min\lbrace p, q\rbrace.
\end{equation*}
The limiting process of this GARCH process is first derived in \cite{BerkesHorvathKokoszka2005} by imposing the assumption $\kappa \in (1/2,1)$. Extending their results, we obtain the limiting process that applies to parameter values that covers the whole range of $(0,1)$. This is a non-trivial extension because when the process deviates further from the integrated GARCH process, the approximation errors in \cite{BerkesHorvathKokoszka2005} diverges and thus a different normalization is needed. 

\section{Main Results}
The main results are summarized in the following one proposition and three theorems. The first proposition modifies the additive representation for $\sigma_t^2$ in \cite{BerkesHorvathKokoszka2005} to accommodate $\kappa \in (0,1)$. Based on the proposition, we establish three theorems to describe the asymptotic behaviours of $\sigma_t^2$ and $u_t$ under the cases $\gamma_n \lesseqqgtr 0$ respectively.

To establish the additive representation of $\sigma_t^2$, we make the following assumptions on the distribution of the innovations $\{\varepsilon_t\}_{t=0}^n$ and the convergence rate of the GARCH coefficients, $\alpha_n$ and $\beta_n$. 
\begin{assumption} \label{ass01}
	$\lbrace \varepsilon_t\rbrace_{t=0}^n$ is an i.i.d sequence with $E\varepsilon_0^2 = 1$ and $E|\varepsilon_0|^{4+\delta} < \infty$, for some $\delta > 0$.
\end{assumption}
\begin{assumption} \label{ass02}
	$\alpha_n\log\log n \rightarrow 0$, $n\alpha_n \rightarrow \infty$ and $\beta_n \rightarrow 1$. 
\end{assumption}
Assumption \ref{ass01} imposes a non-degeneracy condition on the distribution of $\varepsilon_t^2$ and thus ensures its applicability to the central limit theorem. Assumption \ref{ass02} bounds the convergence rate of $\alpha_n$ so that the normalized sequence could converge to a proper limit. Based on these assumptions, we obtain a modified additive representation for $\sigma_t^2$ in Proposition \ref{propAddRep} on the top of \cite{BerkesHorvathKokoszka2005}.
\begin{proposition} [Additive Representation] \label{propAddRep}
	Under Assumption \ref{ass01} and \ref{ass02}, we have the additive representation for $\sigma_t^2$ as
	\begin{align*}
	\sigma_t^2 &= \sigma_0^2t^{t/2}e^{\sqrt{t}\gamma_n}\left(1+\dfrac{\alpha_n}{\sqrt{t}}\sum_{j=1}^{t}\xi_{t-j} + R_{t}^{(1)}\right) + \omega\left[1+\sum_{j=1}^{t}t^{j/2}e^{\frac{j\gamma_n}{\sqrt{t}}}\left(1+\dfrac{\alpha_n}{\sqrt{t}}\sum_{i=1}^{j}\xi_{t-i} + R_{t,j}^{(2)}\right)\left(1+ R_{t,j}^{(3)}\right)\right]
	\end{align*}
	where $\xi_t = \varepsilon_t^{2}-1$ and the remainder terms satisfy
	\begin{align*}
	\left\vert R_{t}^{(1)} \right\vert = O_p\left(\alpha_n^2 + \gamma_n^{2}\right),  &\quad \max\limits_{1 \leq j \leq t}\left\vert R_{t,j}^{(2)} \right\vert = O_p\left(\alpha_n^2\right) \\
	\max\limits_{1 \leq j \leq t}\dfrac{1}{j\log\log j}\left\vert R_{t,j}^{(2)} \right\vert = O_p\left(\dfrac{\alpha_n^2}{t}\right), &\quad \max\limits_{1 \leq j \leq t}\dfrac{1}{j}\left\vert R_{t,j}^{(3)} \right\vert = O_p\left(\dfrac{\alpha_n^2 + \gamma_n^{2}}{t}\right)
	\end{align*}
\end{proposition}

\begin{remark}
	The key difference between our results and \cite{BerkesHorvathKokoszka2005} is the convergence rate of the approximation errors. In \cite{BerkesHorvathKokoszka2005}, the approximation error $\vert R_{t}^{(p)} \vert$, $\forall p = \{1, 2, 3\}$ is of order $t(\alpha_n^2 + \gamma_n^{2})$ or $t\alpha_n^2$ asymptotically. Hence, these errors are negligible only when $\kappa \in (1/2,1)$. We relax this restrictive assumption by normalizing the original terms with $\sqrt{t}$. Under this new normalization, all the approximation errors remains negligible when $\kappa \in (0,1)$.
\end{remark}

To formulate the theorems below, I introduce the following notations. For $0 < t_1 < t_2 < \cdots < t_N <1$ define $k(m) = \lfloor nt_m \rfloor$, $1 \leq m \leq N$. Further, we need the assumptions for relative convergence rate between $\alpha_n$ and $\gamma_n$ to regulate the asymptotic behaviours of returns and volatilities for near-stationary case.

\begin{assumption} \label{ass03}
	$\dfrac{\sqrt{\lvert \gamma_n \rvert}}{\alpha_n n^{1/4}} \rightarrow \infty$, while $\dfrac{\sqrt{\lvert \gamma_n \rvert^{3}}}{\alpha_n n^{1/4}} \rightarrow 0$, as $n \rightarrow \infty$.
\end{assumption}

Assumption \ref{ass03} imposes a rate condition on the localized parameters $\alpha_n$ and $\gamma_n$. This condition is less restrictive than that in \cite{BerkesHorvathKokoszka2005} in the sense that instead of requiring $\lvert \gamma_n \rvert^{3/2}/\alpha_n$ to converge to 0, we allow it to diverge slowly at a rate of $n^{1/4}$. The relaxation of the assumption also attributes to the change of the normalization. 

\begin{theorem} [Near-stationary Case] \label{thmNS}
	Suppose $\gamma_n < 0$, then under Assumption \ref{ass01}-\ref{ass03}, the random variables
	\begin{equation*}
	\dfrac{\sqrt{2\lvert \gamma_n \rvert^{3}}}{\alpha_nk(m)^{1/4}} \dfrac{1}{\sqrt{E\xi_0^2}} \left(\dfrac{\sigma_{k(m)}^2}{\omega {k(m)}^{k(m)/2}} -  \sum_{j=1}^{k(m)-1}e^{\frac{j\gamma_n}{\sqrt{k(m)}}}\right) \xrightarrow{d} \mathcal{N}(0,1).
	\end{equation*}
	
	In addition, the random variables
	\begin{equation*}
	\left(\dfrac{\lvert\gamma_n\rvert }{\omega {k(m)}^{(k(m)+1)/2}} \right)^{1/2} u_{k(m)}
	\end{equation*}
	are asymptotically independent, each with the asymptotic distribution equals to that of $\varepsilon_0$.
\end{theorem}

\begin{theorem} [Integrate Case] \label{thmIN}
	Suppose $\gamma_n = 0$, then under Assumption \ref{ass01} and \ref{ass02}, the volatility has the  asymptotic distribution
	\begin{equation*}
	\dfrac{k(m)^{1/2}}{n^{3/2}\alpha_n}\dfrac{1}{\sqrt{E\xi_0^{2}}}\left(\dfrac{\sigma_{k(m)}^2}{\omega {k(m)}^{k(m)/2}} - k(m)\right)  \xrightarrow{d} \int_{0}^{t_m}xdW(x).
	\end{equation*}
	
	In addition, the random variables
	\begin{equation*}
	\left(\omega {k(m)}^{k(m)/2 + 1}\right)^{-1/2}u_{k(m)}
	\end{equation*}
	are asymptotically independent, each with the asymptotic distribution equals to that of $\varepsilon_0$.
\end{theorem}

Similar to the near-stationary case, we have to impose additional assumption on the relative speed of converging to zero between $\alpha_n$ and $\gamma_n$. 
\begin{assumption} \label{ass04}
	$ \gamma_n / \alpha_n \rightarrow 0$, as $n \rightarrow \infty$.
\end{assumption}
\begin{theorem} [Near-explosive Case] \label{thmNE}
	Suppose $\gamma_n > 0$, then under Assumption \ref{ass01}, \ref{ass02} and \ref{ass04}, the volatility has the asymptotic distribution
	\begin{equation*}
	\dfrac{\gamma_n e^{-\sqrt{k(m)}\gamma_n}}{\alpha_n\sqrt{k(m)}} \dfrac{1}{\sqrt{E\xi_0^2}} \left(\dfrac{\sigma_{k(m)}^2}{\omega k(m)^{k(m)/2}} - \sum_{j=1}^{k(m)-1}e^{\frac{j\gamma_n}{\sqrt{k(m)}}}\right) \Rightarrow W(t_m).
	\end{equation*}
	
	In addition, the random variables
	\begin{equation*}
	\left(\dfrac{\gamma_n e^{-\sqrt{k(m)}\gamma_n}}{\omega k(m)^{(k(m)+1)/2}}\right)^{1/2}u_{k(m)}
	\end{equation*}
	are asymptotically independent, each with the asymptotic distribution equals to that of $\varepsilon_0$.
\end{theorem}

\begin{remark}
	As one may notice, the rate of convergence for both volatility process and return process in all three cases decreases to 0 asymptotically. These seemingly awkward results are reasonable in the sense that the convergence rate is a part of the normalization which reflects the order of the process. In other words, when we compute a partial sum of $X$s in form of $\sum_{i=1}^{n} a_i X_i$, the normalization just plays the role of $a_i$ which is usually required to decrease to 0 for applying a central limit theorem.
\end{remark}

\section{Proofs}
In this section, I present detailed proofs for all the propositions and the theorems listed in the previous section. For readers' convenience, I provide a roadmap for understanding the proofs of the theorems. In general, the proofs are done in three steps:

\textit{Step 1:} We decompose the volatility process into 4 components, $\sigma_{k,s}^2$, $s = 1, \cdots, 4$, by expanding the multiplicative form provided in Proposition \ref{propAddRep}.

\textit{Step 2:} We show the first 3 volatility components are negligible after normalization, and the last term converges to a proper limit by using Cramer-Wold device and Liapounov central limit theorem or Donsker's theorem.

\textit{Step 3:} We figure out a normalization to make the normalized volatility converges to 1. Then, applying this normalization to the return process, we complete the proof.

\begin{proof}[\textbf{Proof of Proposition \ref{propAddRep}}]
	First, note the GARCH(1,1) model can be written into the following multiplicative form:
	\begin{align*}
		\sigma_t^2 &= \sigma_0^2\prod_{i=1}^{t}\left(\beta_n + \alpha_n\varepsilon_{t-i}^2\right) + \omega\left[1 + \sum_{j=1}^{t-1}\prod_{i=1}^{j}\left(\beta_n + \alpha_n\varepsilon_{t-i}^{2}\right)\right] \\
		&= \sigma_0^2 t^{t/2} \prod_{i=1}^{t}\dfrac{\left(\beta_n + \alpha_n\varepsilon_{t-i}^2\right)}{\sqrt{t}} + \omega\left[1 + t^{t/2}\sum_{j=1}^{t-1}\prod_{i=1}^{j}\dfrac{\left(\beta_n + \alpha_n\varepsilon_{t-i}^{2}\right)}{\sqrt{t}}\right].
	\end{align*}
	
	Note that 
	\begin{equation*}
		\max\limits_{1 \leq i \leq t} \dfrac{\left\vert \beta_n + \alpha_n\varepsilon_{t-i}^2 - 1 \right\vert}{\sqrt{t}} \leq \dfrac{\vert \gamma_n \vert}{\sqrt{t}} + \alpha_n\max\limits_{1 \leq i \leq t}\dfrac{\vert \varepsilon_{t-i}^2 - 1 \vert}{\sqrt{t}} =  \dfrac{\vert \gamma_n \vert}{\sqrt{t}} + \alpha_n\max\limits_{1 \leq i \leq t-1}\dfrac{\vert \varepsilon_{i}^2 - 1 \vert}{\sqrt{t}}.
	\end{equation*}
	Then by Assumption \ref{ass01} and \cite{ChowTeicher2012}, we have the almost sure convergence of
	\begin{equation*}
		\max_{1 \leq j \leq t-1}\vert \varepsilon_{i}^2 - 1 \vert = O(\sqrt{t}).
	\end{equation*}
	Therefore, the term above is
	\begin{equation*}
		\max\limits_{1 \leq i \leq t} \dfrac{\left\vert \beta_n + \alpha_n\varepsilon_{t-i}^2 - 1 \right\vert}{\sqrt{t}} = o_p(1).
	\end{equation*}
	
	Now consider the sequence of events 
	$$ A_n = \left\lbrace \max\limits_{1 \leq i \leq t} \dfrac{\lvert\beta_n + \alpha_n\varepsilon_{t-i}^2 - 1 \rvert}{\sqrt{t}} \leq \dfrac{1}{2}\right\rbrace. $$
	From the previous result we know $\lim\limits_{n \rightarrow \infty}P(A_n) = 1$. Then by Taylor expansion, $\lvert \log(1+x) - x \rvert \leq 2x^2$, $\lvert x\rvert \leq 1/2$ on the event $A_n$, which implies
	\begin{align*}
		\left\vert R_{t,j}^{(3)} \right\vert &= \left\vert \sum_{i=1}^{j} \log\dfrac{\left(\beta_n + \alpha_n\varepsilon_{t-i}^2\right)}{\sqrt{t}} - \sum_{i=1}^{j}\dfrac{\left(\gamma_n + \alpha_n\xi_{t-i}\right)}{\sqrt{t}} \right\vert \\
		&= \left\vert \sum_{i=1}^{j} \log\dfrac{\left(\gamma_n + \alpha_n\xi_{t-i} + 1\right)}{\sqrt{t}} - \sum_{i=1}^{j}\dfrac{\left(\gamma_n + \alpha_n\xi_{t-i}\right)}{\sqrt{t}} \right\vert \\
		&\leq \sum_{i=1}^{j} \left\vert \log\left(\dfrac{\gamma_n + \alpha_n\xi_{t-i}}{\sqrt{t}} + 1 \right)- \dfrac{\left(\gamma_n + \alpha_n\xi_{t-i}\right)}{\sqrt{t}} \right\vert \\
		&\leq 2 \sum_{i=1}^{j} \dfrac{\left(\gamma_n + \alpha_n\xi_{t-i}\right)^2}{t} \leq \dfrac{4j\gamma_n^2}{t} + \dfrac{4\alpha_n^2\sum_{i=1}^{j}\xi_{t-i}^2}{t}.
	\end{align*}
	
	By Assumption \ref{ass01} and law of large numbers (LLN), we know
	\begin{equation*}
		\max\limits_{1 \leq j \leq t} \dfrac{1}{j} \left\vert \sum_{i=1}^{j} \xi_{t-i}^2\right\vert \sim \max\limits_{1 \leq j \leq t} \dfrac{1}{j} \left\vert \sum_{i=1}^{j} \xi_{i}^2\right\vert = O_p(1).
	\end{equation*}
	Then by the equation above, we have
	\begin{equation*}
		\max\limits_{1 \leq i \leq j} \dfrac{1}{j}\lvert R_{t,j}^{(3)} \rvert = O_p\left(\dfrac{\gamma_n^2 + \alpha_n^2}{t}\right).
	\end{equation*}
	
	Now by direct plugging into the key multiplicative term we care about, we have 
	\begin{align*}
		\prod_{i=1}^{j}\dfrac{\left(\beta_n + \alpha_n\varepsilon_{t-i}^2\right)}{\sqrt{t}} &= \exp\left\lbrace \sum_{i=1}^{j} \log\left(\dfrac{\beta_n + \alpha_n\varepsilon_{t-i}^2}{\sqrt{t}}\right)\right\rbrace \\
		&= \exp\left\lbrace \dfrac{j\gamma_n}{\sqrt{t}}\right\rbrace \exp\left\lbrace\dfrac{\alpha_n\sum_{i=1}^{j}\xi_{t-i}}{\sqrt{t}} \right\rbrace \exp\left\lbrace R_{t,j}^{(3)}\right\rbrace \\
		&= e^{\frac{j\gamma_n}{\sqrt{t}}}\exp\left\lbrace\dfrac{\alpha_n\sum_{i=1}^{j}\xi_{t-i}}{\sqrt{t}} \right\rbrace\left(1+ R_{t,j}^{(3)}\right).
	\end{align*}
	
	Further, note $\lbrace \xi_t\rbrace_{t=1}^{n}$ is an i.i.d sequence with $E\xi_0^{2} < \infty$, then we know
	\begin{equation*}
		\max\limits_{1 \leq j \leq t} \left\vert \sum_{i=1}^{j} \xi_{t-i} \right\vert = O_p(\sqrt{t}),
	\end{equation*}
	which implies
	\begin{equation*}
		\max\limits_{1 \leq j \leq t} \left\vert \dfrac{\alpha_n}{\sqrt{t}}\sum_{i=1}^{j} \xi_{t-i} \right\vert = O_p(\alpha_n) = o_p(1).
	\end{equation*}
	
	Similarly, we define the sequence of events
	\begin{equation*}
		B_n = \left\lbrace \max\limits_{1 \leq j \leq t}\left\vert \dfrac{\alpha_n}{\sqrt{t}}\sum_{i=1}^{j} \xi_{t-i} \right\vert \leq \dfrac{1}{2} \right\rbrace,
	\end{equation*}
	which is known to have the property $\lim\limits_{n \rightarrow \infty} P(B_n) = 1$. Then by Taylor expansion, $\lvert \exp(x) - (1+x) \rvert \leq \sqrt{e}x^2/2 $ when $\lvert x \rvert \leq 1/2$, on the event $B_n$
	\begin{equation*}
		\left\vert R_{t,j}^{(2)}\right\vert = \left\vert\exp\left\lbrace \dfrac{\alpha_n}{\sqrt{t}}\sum_{i=1}^{j} \xi_{t-i} \right\rbrace - \left(1 + \dfrac{\alpha_n}{\sqrt{t}}\sum_{i=1}^{j} \xi_{t-i}\right)\right\vert \leq \dfrac{\sqrt{e}}{2}\left(\dfrac{\alpha_n}{\sqrt{t}}\sum_{i=1}^{j} \xi_{t-i}\right)^2 = O_p\left(\alpha_n^2\right),
	\end{equation*}
	and by law of iterated logarithm, we know
	\begin{equation*}
		\max\limits_{1 \leq j \leq t}\dfrac{1}{j \log\log j}\left(\dfrac{\alpha_n}{\sqrt{t}}\sum_{i=1}^{j} \xi_{t-i}\right)^2 = O_p\left(\dfrac{\alpha_n^2}{t}\right).
	\end{equation*}
	
	Combining the results above, we have thus showed that
	\begin{equation*}
		\prod_{i=1}^{j}\left(\dfrac{\beta_n + \alpha_n\varepsilon_{t-i}^2}{\sqrt{t}}\right) = e^{\frac{j\gamma_n}{\sqrt{t}}}\left(1 + \dfrac{\alpha_n}{\sqrt{t}}\sum_{i=1}^{j}\xi_{t-i} + R_{t,j}^{(2)}\right)\left(1 + R_{t,j}^{(3)}\right).
	\end{equation*}
	
	Lastly, by the equation above, we know
	\begin{align*}
		\prod_{i=1}^{t}\left(\dfrac{\beta_n + \alpha_n\varepsilon_{t-i}^2}{\sqrt{t}}\right) &= e^{\frac{t\gamma_n}{\sqrt{t}}}\left(1 + \dfrac{\alpha_n}{\sqrt{t}}\sum_{i=1}^{t}\xi_{t-i} + O_p(\alpha_n^2)\right)\left(1 + O_p(\gamma_n^2 + \alpha_n^2)\right) \\
		&= e^{\sqrt{t}\gamma_n}\left(1 + \dfrac{\alpha_n}{\sqrt{t}}\sum_{i=1}^{t}\xi_{t-i} + O_p(\gamma_n^2 + \alpha_n^2)\right),
	\end{align*}
	and this establishes $R_{t}^{(1)}$.
\end{proof}

\begin{proof}[\textbf{Proof of Theorem \ref{thmNS}}]
	First, we focus on the volatilities. Denote $k = \lfloor nt \rfloor$, $0 < t \leq 1$,
	\begin{align*}
		\sigma_k^2 &= \omega + \sigma_0^2 k^{k/2} e^{\sqrt{k}\gamma_n}\left(1+ \dfrac{\alpha_n}{\sqrt{k}}\sum_{j=1}^{k}\xi_{k-j} + R_k^{(1)}\right)  + \omega k^{k/2}\sum_{j=1}^{k-1}e^{\frac{j\gamma_n}{\sqrt{k}}}\left(1+\dfrac{\alpha_n}{\sqrt{k}}\sum_{i=1}^{j}\xi_{k-i} + R_{k,j}^{(2)}\right)R_{k,j}^{(3)} \\
		&\ \ \ + \omega k^{k/2} \sum_{j=1}^{k-1}e^{\frac{j\gamma_n}{\sqrt{k}}} R_{k,j}^{(2)} + \omega k^{k/2} \sum_{j=1}^{k-1}e^{\frac{j\gamma_n}{\sqrt{k}}}\left(1+\dfrac{\alpha_n}{\sqrt{k}}\sum_{i=1}^{j}\xi_{k-i}\right) \\
		&= \omega + \sigma_{k,1}^{2} + \sigma_{k,2}^{2} + \sigma_{k,3}^{2} + \sigma_{k,4}^{2}.
	\end{align*}
	
	For $\sigma_{k,1}^{2}$, note $k^{-1/2}\sum_{j=1}^{k}\xi_{k-j}$ is asymptotically normal, then by Proposition \ref{propAddRep},
	\begin{equation*}
		\dfrac{\alpha_n}{\sqrt{k}}\sum_{j=1}^{k}\xi_{k-j} + R_{k}^{(1)} = o_p(1),
	\end{equation*}
	and this implies
	\begin{align*}
		\left\lvert \sigma_{k,1}^2 \right\rvert &= O_p\left(k^{k/2}e^{\sqrt{k}\gamma_n}\right).
	\end{align*}
	
	For $\sigma_{k,2}^{2}$, note by Lemma 4.1 in \cite{BerkesHorvathKokoszka2005}, we have
	\begin{equation} \label{eqApp01}
		\sum_{j=1}^{k}je^{\frac{j\gamma_n}{\sqrt{k}}} \sim \dfrac{k}{\lvert \gamma_n \rvert^2}\Gamma(2),
	\end{equation}
	and note that
	\begin{equation} \label{eqApp02}
		\max\limits_{1\leq j \leq k-1} \left\lvert \dfrac{\alpha_n}{\sqrt{k}}\sum_{i=1}^{j}\xi_{k-i} + R_{k,j}^{(2)} \right\rvert = o_p(1).
	\end{equation}
	Then by equation (\ref{eqApp01}), (\ref{eqApp02}) and Proposition \ref{propAddRep} we have
	\begin{align*}
		\left\lvert \sigma_{k,2}^{2} \right\rvert &= \left\lvert \omega k^{k/2}\sum_{j=1}^{k-1}je^{\frac{j\gamma_n}{\sqrt{k}}}\left(1 + \dfrac{\alpha_n}{\sqrt{k}}\sum_{i=1}^{j}\xi_{k-i} + R_{k,j}^{(2)}\right)\dfrac{1}{j}R_{k,j}^{(3)} \right\rvert \\
		&= O_p(1)\omega k^{k/2} \dfrac{\alpha_n^2 + \gamma_n^2}{k} \dfrac{k}{\lvert \gamma_n \rvert^2} \\
		&= O_p\left(\dfrac{k^{k/2}\left(\alpha_n^2 + \gamma_n^2\right)}{\gamma_n^2}\right).
	\end{align*}
	
	For $\sigma_{k,3}^{2}$, similarly, by Proposition \ref{propAddRep} and Lemma 4.1 in \cite{BerkesHorvathKokoszka2005}, we have
	\begin{align*}
		\left\lvert \sigma_{k,3}^{2} \right\rvert &= \left\lvert \omega k^{k/2}\sum_{j=1}^{k-1}e^{\frac{j\gamma_n}{\sqrt{k}}}R_{k,j}^{(2)} \right\rvert \\
		&= O_p(1)\omega k^{k/2}\dfrac{\alpha_n^2}{k}\sum_{j=1}^{k-1}je^{\frac{j\gamma}{\sqrt{k}}}\log\log j \\
		&= O_p\left(\dfrac{k^{k/2} \left(\alpha_n^{2}\log\log k\right) }{\gamma_n^2}\right).
	\end{align*}
	
	Lastly, for $\sigma_{k,4}^{2}$, by Lemma 4.1 in \ref{propAddRep} we have
	\begin{align*}
		\sigma_{k,4}^{2} &= \omega k^{k/2} \sum_{j=1}^{k-1}e^{\frac{j\gamma_n}{\sqrt{k}}} + \omega k^{k/2} \dfrac{\alpha_n}{\sqrt{k}} \sum_{j=1}^{k-1}e^{\frac{j\gamma_n}{\sqrt{k}}} \sum_{i=1}^{j}\xi_{k-i} \\
		&= O_p\left(\dfrac{k^{k/2}k^{1/2}}{\lvert \gamma_n \rvert}\right) + \omega k^{k/2} \dfrac{\alpha_n}{\sqrt{k}} \sum_{j=1}^{k-1}e^{\frac{j\gamma_n}{\sqrt{k}}} \sum_{i=1}^{j}\xi_{k-i}.
	\end{align*}
	
	Therefore, we only have to consider the last term in the above equation. Define 
	\begin{equation*}
		\tau_m = k(m)^{-1/4}\sum_{j=1}^{k(m)-1}e^{\frac{j\gamma_n}{\sqrt{k(m)}}} \xi_{k(m)-j}, \quad 1 \leq m \leq N,
	\end{equation*}
	and
	\begin{equation*}
		\tau_m^* = k(m)^{-1/2}\sum_{j=1}^{k(m)-1}e^{\frac{j\gamma_n}{\sqrt{k(m)}}} \sum_{i=1}^{j}\xi_{k(m)-i}, \quad 1 \leq m \leq N.
	\end{equation*}
	
	Then by Cramer-Wold device (Theorem 29.4 of \cite{Billingsley1995}), we have
	\begin{align*}
		\sum_{m=1}^{N}\mu_m\tau_m &= \sum_{i=1}^{k(1)-1}\sum_{m=1}^{N}\dfrac{\mu_m}{k(m)^{1/4}} e^{\frac{(k(m)-i)\gamma_n}{\sqrt{k(m)}}} + \sum_{i=k(1)}^{k(2)-1}\sum_{m=2}^{N}\dfrac{\mu_m}{k(m)^{1/4}} e^{\frac{(k(m)-i)\gamma_n}{\sqrt{k(m)}}}  \\
		&+ \cdots + \sum_{i=k(N-1)}^{k(N)-1}\dfrac{\mu_N}{k(N)^{1/4}} e^{\frac{(k(N)-i)\gamma_n}{\sqrt{k(N)}}} \\
		&= S_1 + S_2 + \cdots + S_N.
	\end{align*}
	
	Observe that
	\begin{align*}
		ES_1^2 &= E\xi_0^2\left( \sum_{i=1}^{k(1)-1}\sum_{m=1}^{N}k(m)^{-1/4}\mu_m e^{\frac{(k(m)-i)\gamma_n}{\sqrt{k(m)}}}\right)^2 \\
		&= E\xi_0^2 \sum_{m=1}^{N}\dfrac{\mu_m^2}{\sqrt{k(m)}}\sum_{i=1}^{k(1)-1} e^{\frac{2(k(m)-i)\gamma_n}{\sqrt{k(m)}}} + E\xi_0^2 \sum_{1 \leq m \neq l \leq N}\left(k(m)k(l)\right)^{-1/4}\mu_m\mu_l\sum_{i=1}^{k(1)-1} e^{\frac{(k(m)-i)\gamma_n}{\sqrt{k(m)}}+\frac{(k(l)-i)\gamma_n}{\sqrt{k(l)}}} \\
		&= E\xi_0^2 \dfrac{\mu_1^2}{\sqrt{k(1)}}\sum_{i=1}^{k(1)-1} e^{\frac{2(k(1)-i)\gamma_n}{\sqrt{k(1)}}} + E\xi_0^2 \sum_{m=2}^{N}\dfrac{\mu_m^2}{\sqrt{k(m)}}\sum_{i=1}^{k(1)-1} e^{\frac{2(k(m)-i)\gamma_n}{\sqrt{k(m)}}} \\
		&\ \ \ + E\xi_0^2 \sum_{1 \leq m \neq l \leq N}\left(k(m)k(l)\right)^{-1/4}\mu_m\mu_l\sum_{i=1}^{k(1)-1} e^{\frac{(k(m)-i)\gamma_n}{\sqrt{k(m)}}+\frac{(k(l)-i)\gamma_n}{\sqrt{k(l)}}} \\
		&= E\xi_0^2 \dfrac{\mu_1^2}{\sqrt{k(1)}}\sum_{i=1}^{k(1)-1} e^{\frac{2i\gamma_n}{\sqrt{k(1)}}} + E\xi_0^2 \sum_{m=2}^{N}\dfrac{\mu_m^2}{\sqrt{k(m)}}e^{\frac{2(k(m)-k(1))\gamma_n}{\sqrt{k(m)}}}\sum_{i=1}^{k(1)-1} e^{\frac{2i\gamma_n}{\sqrt{k(m)}}} \\
		&\ \ \ + E\xi_0^2 \sum_{1 \leq m \neq l \leq N}(k(m)k(l))^{-1/4}\mu_m\mu_le^{\frac{(k(m)-k(1))\gamma_n}{\sqrt{k(m)}}+\frac{(k(l)-k(1))\gamma_n}{\sqrt{k(l)}}}\sum_{i=1}^{k(1)-1} e^{\frac{i\gamma_n}{\sqrt{k(m)}}+\frac{i\gamma_n}{\sqrt{k(l)}}} \\
		&\sim E\xi_0^2 \mu_1^2\dfrac{1}{2\lvert \gamma_n\rvert} + E\xi_0^2 \sum_{m=2}^{N}\mu_m^2e^{\frac{2(k(m)-k(1))\gamma_n}{\sqrt{k(m)}}} \dfrac{1}{2\lvert \gamma_n \rvert} \\
		&\ \ \ + E\xi_0^2 \sum_{1 \leq m \neq l \leq N}\dfrac{\mu_m\mu_l}{(\sqrt{k(m)} +\sqrt{k(l)}) \lvert\gamma_n\rvert} e^{\frac{(k(m)-k(1))\gamma_n}{\sqrt{k(m)}} + \frac{(k(l)-k(1))\gamma_n}{\sqrt{k(l)}}}  \\
		&= E\xi_0^2 \mu_1^2\dfrac{1}{2\lvert \gamma_n\rvert} + o\left(\dfrac{1}{\lvert \gamma_n \rvert}\right),
	\end{align*}
	we then have
	\begin{align*}
		E\left(\sum_{m=1}^{N} \mu_m \tau_m\right)^2 &= \left(\sum_{m=1}^{N} \mu_m^2\right) E\xi_0\dfrac{1}{2\lvert \gamma_n \rvert} + o\left(\dfrac{1}{\lvert \gamma_n \rvert}\right).
	\end{align*}
	
	Observe also that, for some $c_i$, $1 \leq i \leq k(N)-1$, we have
	\begin{equation*}
		\sum_{m=1}^{N} \mu_m \tau_m = \sum_{i=1}^{k(N)-1} c_i \xi_i,
	\end{equation*}
	and by Jensen's inequality, we know for some $\delta > 0$,
	\begin{align*}
		\lvert c_i \rvert^{2+\delta} &= \left\vert k(1)^{-1/4}\mu_1 e^{\frac{(k(1)-i) \gamma_n}{\sqrt{k(1)}}} + k(2)^{-1/4}\mu_1 e^{\frac{(k(2)-i) \gamma_n}{\sqrt{k(2)}}}+ \cdots +k(N)^{-1/4}\mu_1 e^{\frac{(k(N)-i) \gamma_n}{\sqrt{k(N)}}} \right\vert^{2+\delta} \\
		&\leq C_1(N)\left[\dfrac{\lvert \mu_1\rvert^{2+\delta}}{k(1)^{1/2+\delta/4}}e^{\frac{(k(1) - i)(2+\delta)\gamma_n}{\sqrt{k(1)}}} + \cdots + \dfrac{\lvert \mu_N\rvert^{2+\delta}}{k(N)^{1/2+\delta/4}}e^{\frac{(k(N) - i)(2+\delta)\gamma_n}{\sqrt{k(N)}}}\right].
	\end{align*}
	This implies that
	\begin{equation*}
		\sum_{i=1}^{k(N)-1} \lvert c_i \rvert^{2+\delta} \sim C_1(N)\lvert \mu_1\rvert^{2+\delta} \dfrac{1}{k(1)^{\delta/4}(2+\delta)\lvert \gamma_n \rvert} + O\left(\dfrac{1}{k(2)^{\delta/4}\lvert \gamma_n\rvert}\right) = o\left(\dfrac{1}{\lvert \gamma_n \rvert}\right).
	\end{equation*}
	
	Now we can easily check the Liapounov's condition, where
	\begin{equation*}
		\dfrac{\left(\sum_{i=1}^{k(N)-1}\lvert c_i \rvert^{2+\delta} E\lvert \xi_i\rvert^{2+\delta}\right)^{1/(2+\delta)}}{\left(\sum_{i=1}^{k(N)-1}c_i^2 E\xi_i^{2}\right)^{1/2}} = o\left(\lvert \gamma_n \rvert^{1/2 - 1/(2+\delta)}\right) = o_p(1).
	\end{equation*}
	Then by Liapounov central limit theorem (Theorem 27.3, p.362 of \cite{Billingsley1995}), we have
	\begin{equation*}
		\sqrt{2\lvert \gamma_n \rvert}\left[\tau_1, \tau_2, \cdots, \tau_N\right] \xrightarrow{d} \sqrt{E\xi_{0}^2}\left[\eta_1, \eta_2, \cdots, \eta_N\right],
	\end{equation*}
	where $\eta_1, \eta_2, \cdots, \eta_N$ are independent standard normal random variables.
	
	Now we have to check the relationship between $\tau_m$ and $\tau_m^*$. Note by $k^{-1/2}\left(e^{\frac{\gamma_n}{\sqrt{k}}}-1\right)^{-1} = \left(\gamma_n + o(1)\right)^{-1}$, we have
	\begin{align*}
		\dfrac{1}{\sqrt{k}}\sum_{j=i}^{k-1}e^{\frac{j\gamma_n}{\sqrt{k}}} - \lvert \gamma_n\rvert^{-1}e^{\frac{i\gamma_n}{\sqrt{k}}} &= \dfrac{1}{\sqrt{k}}\dfrac{e^{\frac{k\gamma_n}{\sqrt{k}}} - e^{\frac{i\gamma_n}{\sqrt{k}}}}{e^{\frac{\gamma_n}{\sqrt{k}}}-1} - \lvert \gamma_n\rvert^{-1}e^{\frac{i\gamma_n}{\sqrt{k}}} \\
		&= \left(\gamma_n + o(1)\right)^{-1}\left(e^{\frac{k\gamma_n}{\sqrt{k}}} - e^{\frac{i\gamma_n}{\sqrt{k}}}\right) - \lvert \gamma_n \rvert^{-1} e^{\frac{i\gamma_n}{\sqrt{k}}} \\
		&= \left(\gamma_n^{-1} + O(1)\right) e^{\frac{k\gamma_n}{\sqrt{k}}} - e^{\frac{i\gamma_n}{\sqrt{k}}}O(1).
	\end{align*}
	
	Then, we know
	\begin{align*}
		E\left[\sqrt{2\lvert \gamma_n\rvert^3}\tau_m^* - \sqrt{2\lvert \gamma_n\rvert}\tau_m\right]^2 &= \dfrac{2\lvert \gamma_n \rvert^3}{\sqrt{k}} E\left[\dfrac{1}{\sqrt{k}}\sum_{i=1}^{k-1}\left(\sum_{j=i}^{k-1}e^{\frac{j\gamma}{\sqrt{k}}}\right)\xi_{k-i} - \lvert \gamma_n\rvert^{-1}\sum_{i=1}^{k-1}e^{\frac{i\gamma_n}{\sqrt{k}}}\xi_{k-i} \right]^{2} \\
		&= \dfrac{2\lvert \gamma_n \rvert^3}{\sqrt{k}} E\xi_0^2\sum_{i=1}^{k-1}\left(\dfrac{1}{\sqrt{k}}\sum_{j=i}^{k-1}e^{\frac{j\gamma_n}{\sqrt{k}}} - \lvert \gamma_n\rvert^{-1}e^{\frac{i\gamma_n}{\sqrt{k}}}\right)^2 \\
		&\sim \dfrac{2\lvert \gamma_n \rvert^3}{\sqrt{k}} E\xi_0^2\left(k\gamma_n^{-2}e^{\sqrt{k}\gamma_n} + \dfrac{\sqrt{k}}{2\lvert\gamma\rvert} - 2\gamma_n^{-1}e^{\sqrt{k}\gamma_n}\dfrac{\sqrt{k}}{\lvert\gamma\rvert}\right) \\
		&= 2E\xi_0^2 O\left(\sqrt{k}\lvert\gamma_n\rvert e^{\sqrt{k}\gamma_n}\right) + o_p(1) \\
		&= o_p(1),
	\end{align*}
	where the last equality comes from the well known limits of $xe^{-x}$,
	\begin{equation*}
		\lim\limits_{x \rightarrow \infty} \dfrac{x}{e^{x}} = \lim\limits_{x \rightarrow \infty}\dfrac{1}{e^x} = 0 \quad \text{and} \quad \lim\limits_{x \rightarrow 0} \dfrac{x}{e^{x}} = 0.
	\end{equation*}
	
	Therefore, we have
	\begin{equation*}
		\sqrt{2\lvert \gamma_n \rvert^{3}}\left[\tau_1^*, \tau_2^*, \cdots, \tau_N^* \right] \xrightarrow{d} \sqrt{E\xi_0^2}\left[\eta_1, \eta_2, \cdots, \eta_N\right],
	\end{equation*}
	
	Now combine the results above, we have, for each $k = {\lfloor nt_m \rfloor}$, $m = 1, \cdots, N$
	\begin{align*}
		\dfrac{\sqrt{2\lvert \gamma_n \rvert^{3}}}{\alpha_n k^{1/4}} \dfrac{1}{\sqrt{E\xi_0^2}} \left(\dfrac{\sigma_{k}^2}{\omega {k}^{k/2}} - \sum_{j=1}^{k-1}e^{\frac{j\gamma_n}{\sqrt{k}}}\right) \xrightarrow{d} \mathcal{N}(0,1).
	\end{align*}
	
	Now, for returns, we know from the above result that
	\begin{align*}
		\dfrac{\lvert\gamma_n\rvert\sigma_{k}^2}{\omega {k}^{(k+1)/2}} - 1 = O_p\left(\dfrac{\alpha_n n^{1/4}}{\sqrt{\lvert \gamma_n \rvert}}\right) = o_p(1).
	\end{align*}
	Therefore, by the return equation, we have
	\begin{equation*}
		\left(\dfrac{\lvert\gamma_n\rvert }{\omega {k}^{(k+1)/2}} \right)^{1/2} u_k =  \left(\dfrac{\lvert\gamma_n\rvert\sigma_{k}^2}{\omega {k}^{(k+1)/2}}\right)^{1/2} \varepsilon_k \sim \varepsilon_k.
	\end{equation*}
\end{proof}

\begin{proof}[\textbf{Proof of Theorem \ref{thmIN}}]
	Similar to Theorem \ref{thmNS}, when $\gamma_n = 0$, the volatility admits the decomposition. Then, for $\sigma_{k,1}^{2}$, by central limit theorem, we know
	\begin{equation*}
		\dfrac{\alpha_n}{\sqrt{k}}\sum_{j=1}^{k}\xi_{k-j} = O_p(\alpha_n) = o_p(1)
	\end{equation*}
	which, combining with Proposition \ref{propAddRep}, implies that
	\begin{equation*}
		\left\lvert \sigma_{k,1}^2 \right\rvert = O_p\left(k^{k/2}\right)
	\end{equation*}
	
	For $\sigma_{k,2}^{2}$, note that we have established equation (\ref{eqApp02}), then by Proposition \ref{propAddRep}, we have
	\begin{equation*}
		\left\lvert \sigma_{k,2}^{2} \right\rvert = O_p\left(k^{k/2}\alpha_n^{2}\right).
	\end{equation*}
	
	For $\sigma_{k,3}^{2}$, by Proposition \ref{propAddRep} we have
	\begin{equation*}
		\left\lvert \sigma_{k,3}^{2} \right\rvert = O_p(k^{k/2}\alpha_n^{2}).
	\end{equation*}
	
	Lastly, for $\sigma_{k,4}^2$, note by Lemma 5.1 in \cite{BerkesHorvathKokoszka2005}, for $k = \lfloor nt \rfloor$, $t \in (0, 1)$ , we have
	\begin{equation*}
		\dfrac{1}{n^{3/2}}\sum_{j=1}^{\lfloor nt \rfloor - 1}\sum_{i=1}^{j}\xi_{k-i} \xrightarrow{d} \sqrt{E\xi_{0}^2}\int_{0}^{t}xdW(x),
	\end{equation*}
	where $W(x)$ is a Wiener process. 
	
	Therefore, for $k(m) = \lfloor nt_m \rfloor$, $m = 1, \cdots, N$, we have
	\begin{align*}
		\dfrac{k(m)^{1/2}}{n^{3/2}\alpha_n}\left(\dfrac{\sigma_{k(m)}^2}{\omega {k(m)}^{k(m)/2}} - k(m)\right) = \dfrac{1}{n^{3/2}}\sum_{j=1}^{\lfloor nt_m \rfloor - 1}\sum_{i=1}^{j}\xi_{k(m)-i} + o_p(1) \xrightarrow{d} \sqrt{E\xi_0^{2}}\int_{0}^{t_m}xdW(x).
	\end{align*}
	
	Further, note the results above implies that
	\begin{equation*}
		\dfrac{\sigma_{k(m)}^2}{\omega {k(m)}^{k(m)/2 + 1}} - 1 = O_p\left(\left(\dfrac{n}{k}\right)^{3/2}\alpha_n\right) = o_p(1).
	\end{equation*}
	Hence, by return equation, we obtain
	\begin{equation*}
		\left(\dfrac{1}{\omega {k(m)}^{k(m)/2 + 1}}\right)^{1/2}u_{k(m)} = \left(\dfrac{\sigma_{k(m)}^2}{\omega {k(m)}^{k(m)/2 + 1}}\right)^{1/2}\varepsilon_{k(m)} \xrightarrow{d} \varepsilon_{k(m)}.
	\end{equation*}
\end{proof}

\begin{proof}[\textbf{Proof of Theorem \ref{thmNE}}]
	Similar to proof of Theorem \ref{thmNS}, when $\gamma_n > 0$, the volatility admits the additive representation. For $\sigma_{k,1}^{2}$, similar to that in Theorem \ref{thmNS}, 
	\begin{equation*}
	\left\lvert \sigma_{k,1}^2 \right\rvert = O_p\left(k^{k/2}e^{\sqrt{k}\gamma_n}\right).
	\end{equation*}
	
	For $\sigma_{k,2}^2$, by Proposition \ref{propAddRep} and equation (\ref{eqApp02}), we have the relation
	\begin{align*}
	\left\lvert \sigma_{k,2}^{2} \right\rvert &= \left\lvert \omega k^{k/2}\sum_{j=1}^{k-1} je^{\frac{j\gamma_n}{\sqrt{k}}}(1+o_p(1))\dfrac{1}{j}R_{k,j}^{(3)} \right\rvert \\
	&=O_p(1)\omega k^{k/2}(\alpha_n^{2} + \gamma_n^{2})\dfrac{e^{\frac{k\gamma_n}{\sqrt{k}}} - e^{\frac{\gamma_n}{\sqrt{k}}}}{e^{\frac{\gamma_n}{\sqrt{k}}}-1} \\
	&< O_p\left(k^{k/2}\left(\alpha_n^2 + \gamma_n^2\right) \dfrac{\sqrt{k}e^{\sqrt{k}\gamma_n}}{\gamma_n}\right),
	\end{align*}
	where the last inequality comes from the fact that
	\begin{equation*}
	\dfrac{e^{\frac{k\gamma_n}{\sqrt{k}}} - e^{\frac{\gamma_n}{\sqrt{k}}}}{e^{\frac{\gamma_n}{\sqrt{k}}}-1} < \dfrac{e^{\sqrt{k}\gamma_n}}{\gamma_n/\sqrt{k}}.
	\end{equation*}
	
	For $\sigma_{k,3}^2$, by Proposition \ref{propAddRep}, we have
	\begin{align*}
	\left\lvert \sigma_{k,3}^2 \right\rvert &= \left\lvert \omega k^{k/2} \sum_{j=1}^{k-1}e^{\frac{j\gamma_n}{\sqrt{k}}} \left(j\log\log j\right) \dfrac{1}{j\log\log j} R_{k,j}^{(2)}\right\rvert \\
	&= O_p(1)\omega k^{k/2} \left(k\log\log k\right) \dfrac{\alpha_n^{2}}{k}\dfrac{e^{\frac{k\gamma_n}{\sqrt{k}}} - e^{\frac{\gamma_n}{\sqrt{k}}}}{e^{\frac{\gamma_n}{\sqrt{k}}}-1} \\
	&< O_p\left(k^{k/2}\left(\alpha_n^2\log\log k\right)\dfrac{\sqrt{k}e^{\sqrt{k}\gamma_n}}{\gamma_n} \right).
	\end{align*}
	
	Lastly, for $\sigma_{k,4}^{2}$, we have
	\begin{equation*}
	\sigma_{k,4}^2 = \omega k^{k/2} \sum_{j=1}^{k-1}e^{\frac{j\gamma_n}{\sqrt{k}}} + \omega k^{k/2} \dfrac{\alpha_n}{\sqrt{k}} \sum_{j=1}^{k-1}e^{\frac{j\gamma_n}{\sqrt{k}}} \sum_{i=1}^{j}\xi_{k-i}.
	\end{equation*}
	
	Now, we introduce the following lemma to assist the proof.
	\begin{lemma} \label{lemApp}
		If Assumption \ref{ass01} and \ref{ass02} hold, then
		\begin{equation*}
		\dfrac{\gamma_n^2}{k}e^{-2\sqrt{k}\gamma_n}E\left(\dfrac{1}{\sqrt{k}}\sum_{j=1}^{k-1}e^{\frac{j\gamma_n}{\sqrt{k}}}\sum_{i=1}^{j}\xi_{k-i} - \dfrac{e^{\sqrt{k}\gamma_n}}{\gamma_n}\sum_{i=1}^{k-1}\xi_i\right)^2 \rightarrow 0.
		\end{equation*}
	\end{lemma}
	
	Then by Lemma \ref{lemApp}, we have
	\begin{align*}
	\dfrac{\gamma_n e^{-\sqrt{k}\gamma_n}}{\sqrt{k}\alpha_n} \left(\dfrac{\sigma_{k,4}^2}{\omega k^{k/2}} - \sum_{j=1}^{k-1}e^{\frac{j\gamma_n}{\sqrt{k}}}\right) &= \dfrac{\gamma_n e^{-\sqrt{k}\gamma_n}}{\sqrt{k}}\dfrac{1}{\sqrt{k}}\sum_{j=1}^{k-1}e^{\frac{j\gamma_n}{\sqrt{k}}}\sum_{i=1}^{j}\xi_{k-i} + o_p(1) = \dfrac{1}{\sqrt{k}}\sum_{i=1}^{k-1}\xi_i + o_p(1).
	\end{align*}
	
	Therefore, by Donsker's theorem, we obtain that, for $k(m) = \lfloor nt_m \rfloor$, $t_m \in (0, 1)$ and $m = 1, 2, \cdots, N$,
	\begin{equation*}
	\dfrac{\gamma_n e^{-\sqrt{k(m)}\gamma_n}}{\sqrt{k(m)}\alpha_n} \dfrac{1}{\sqrt{E\xi_0^2}} \left(\dfrac{\sigma_{k(m)}^2}{\omega k(m)^{k(m)/2}} - \sum_{j=1}^{k(m)-1}e^{\frac{j\gamma_n}{\sqrt{k(m)}}}\right) \Rightarrow W(t_m),
	\end{equation*}
	where $W(t)$ is a finite dimensional Wiener process.
	
	Further, note that
	\begin{equation*}
	\dfrac{\gamma_n}{\sqrt{k}}e^{-\sqrt{k}\gamma_n}\left(\sum_{j=1}^{k-1}e^{\frac{j\gamma_n}{\sqrt{k}}} - \dfrac{\sqrt{k}e^{\sqrt{k}\gamma_n}}{\gamma_n}\right) = o(1),
	\end{equation*}
	then by the result above we know
	\begin{equation*}
	\dfrac{\gamma_n e^{-\sqrt{k(m)}\gamma_n}}{\sqrt{k(m)}}  \left(\dfrac{\sigma_{k(m)}^2}{\omega k(m)^{k(m)/2}} - \sum_{j=1}^{k(m)-1}e^{\frac{j\gamma_n}{\sqrt{k(m)}}}\right) = O_p(\alpha_n) = o_p(1).
	\end{equation*}
	
	Hence, by return equantion, we derive
	\begin{equation*}
	\left(\dfrac{\gamma_n e^{-\sqrt{k}\gamma_n}}{\omega k^{(k+1)/2}}\right)^{1/2}u_{k} = \left(\dfrac{\gamma_n e^{-\sqrt{k}\gamma_n}}{\omega k^{(k+1)/2}}\sigma_k^2\right)^{1/2}\varepsilon_{k} \sim \varepsilon_k.
	\end{equation*}
\end{proof}

\textbf{Proof of Lemma \ref{lemApp}.} Note that 
\begin{equation*}
\dfrac{1}{\sqrt{k}}\sum_{j=1}^{k-1}e^{\frac{j\gamma_n}{\sqrt{k}}}\sum_{i=1}^{j}\xi_{k-i} = \dfrac{1}{\sqrt{k}}\sum_{i=1}^{k-1}\left(\sum_{j=i}^{k-1}e^{\frac{j\gamma_n}{\sqrt{k}}}\right)\xi_{k-i} \quad \text{and} \quad  \sum_{i=1}^{k-1}\xi_i = \sum_{i=1}^{k-1}\xi_{k-i},
\end{equation*}
Then,
\begin{align*}
E\left(\dfrac{1}{\sqrt{k}}\sum_{j=1}^{k-1}e^{\frac{j\gamma_n}{\sqrt{k}}}\sum_{i=1}^{j}\xi_{k-i} - \dfrac{e^{\sqrt{k}\gamma_n}}{\gamma_n}\sum_{i=1}^{k-1}\xi_i\right)^2 &= E\xi_0^{2}\sum_{i=1}^{k-1}\left(\dfrac{1}{\sqrt{k}}\sum_{j=i}^{k-1}e^{\frac{j\gamma_n}{\sqrt{k}}} - \dfrac{\sqrt{k}e^{\sqrt{k}\gamma_n}}{\gamma_n} \right)^2 \\
&= \dfrac{E\xi_0^{2}}{k}\sum_{i=1}^{k-1}\left(\dfrac{e^{\frac{k\gamma_n}{\sqrt{k}}} - e^{\frac{i\gamma_n}{\sqrt{k}}}}{e^{\frac{\gamma_n}{\sqrt{k}}}-1} - \dfrac{\sqrt{k}e^{\sqrt{k}\gamma_n}}{\gamma_n} \right)^2.
\end{align*}

Note by Taylor expansion,
\begin{equation*}
\left\lvert \dfrac{e^{\frac{k\gamma_n}{\sqrt{k}}} - e^{\frac{i\gamma_n}{\sqrt{k}}}}{e^{\frac{\gamma_n}{\sqrt{k}}}-1} - \dfrac{\sqrt{k}e^{\sqrt{k}\gamma_n}}{\gamma_n} \right\rvert \leq C_1\left(\dfrac{\sqrt{k}e^{\frac{i\gamma_n}{\sqrt{k}}}}{\gamma_n} + e^{\sqrt{k}\gamma_n}\right),
\end{equation*}
which implies that
\begin{align*}
\sum_{i=1}^{k-1}\left(\dfrac{e^{\frac{k\gamma_n}{\sqrt{k}}} - e^{\frac{i\gamma_n}{\sqrt{k}}}}{e^{\frac{\gamma_n}{\sqrt{k}}}-1} - \dfrac{\sqrt{k}e^{\sqrt{k}\gamma_n}}{\gamma_n}\right)^2 &\leq 2C_1^2\left(\sum_{i=1}^{k-1}\dfrac{ke^{\frac{2i\gamma_n}{\sqrt{k}}}}{\gamma_n^2} + ke^{2\sqrt{k}\gamma_n} \right) \\
&= O(1)\left(\dfrac{k}{\gamma_n^2}\dfrac{e^{\frac{2k\gamma_n}{\sqrt{k}}} - e^{\frac{2\gamma_n}{\sqrt{k}}}}{e^{\frac{2\gamma_n}{\sqrt{k}}} - 1} + ke^{2\sqrt{k}\gamma_n} \right) \\
&= O(1)\left(\dfrac{k}{\gamma_n^3}e^{2\sqrt{k}\gamma_n} + ke^{2\sqrt{k}\gamma_n}\right).
\end{align*}

Now we can see that
\begin{align*}
\dfrac{\gamma_n^2}{k}e^{-2\sqrt{k}\gamma_n}E\left(\dfrac{1}{\sqrt{k}}\sum_{j=1}^{k-1}e^{\frac{j\gamma_n}{\sqrt{k}}}\sum_{i=1}^{j}\xi_{k-i} - \dfrac{e^{\sqrt{k}\gamma_n}}{\gamma_n}\sum_{i=1}^{k-1}\xi_i\right)^2 &= O(1)\dfrac{\gamma_n^2}{k}e^{-2\sqrt{k}\gamma_n}\dfrac{E\xi_0^2}{k}\left(\dfrac{k}{\gamma_n^3}e^{2\sqrt{k}\gamma_n} + ke^{2\sqrt{k}\gamma_n}\right) = o_p(1).
\end{align*}

\nocite{*}
\bibliography{ExGARCHBib}

\end{document}